\documentclass{amsart}
\usepackage{amsmath}
\usepackage{amsfonts}
\usepackage{amssymb}
\usepackage{graphicx}

 \newtheorem{thm}{Theorem}[section]
 \newtheorem{cor}[thm]{Corollary}
 
 \newtheorem{prop}[thm]{Proposition}
 \theoremstyle{definition}
 
 \theoremstyle{remark}
 \newtheorem{rem}[thm]{Remark}
 
 \numberwithin{equation}{section}

\begin{document}

\title[A generalization of the probability that the commutator of two ...]
 {A generalization of the probability that the commutator of two group elements is equal to a given element}


\author[A. M. A. Alghamdi]{Ahmad M.A. Alghamdi}
\address{Department of Mathematical Sciences, Faculty of Applied Sciences\\
Umm Alqura University, P.O. Box 14035, Makkah, 21955, Saudi Arabia}
\email{amghamdi@uqu.edu.sa}


\author[F. G. Russo]{Francesco G. Russo}
\address{IISSS "Axel Munthe"\\
viale Axel Munthe, 80074, Anacapri (Naples), Italy}
\email{francescog.russo@yahoo.com}




\date{\today}
\dedicatory{}

\begin{abstract}
The probability that the commutator of two group elements is equal to a given element has been introduced in
literature few years ago. Several authors have investigated this notion with  methods of the representation
theory and with combinatorial techniques. Here we illustrate that a wider context may be considered and show
some structural restrictions on the group.
\end{abstract}

\maketitle

\section{Different formulations of the commutativity degree}
Given two elements $x$ and $y$ of a group $G$, several authors
studied the probability that a randomly chosen commutator $[x,y]$ of
$G$ satisfies a prescribed property. P. Erd\H{o}s and P. Tur\'an
\cite{et} began to investigate the case $[x,y]=1$,  noting some
structural restrictions on $G$ from bounds of statistical nature.
Their approach involved combinatorial techniques, which were
developed successively in
 \cite{cms, das-nath2, das-nath1, doostie-maghasedi, elr, err, erovenko-sury, gr, gustafson, l, rezaei1} and
 extended to the infinite case in \cite{er, gustafson, rezaei2}.  On another hand, P. X. Gallagher
\cite{gallagher} investigated the case $[x,y]=1$, using character theory, and opened another line of research,
illustrated in \cite{das-nath2, das-nath1, gr, ps, rusin}. The literature shows that it is possible to variate
the condition on $[x,y]$ involving arbitrary words, which could not be the commutator word $[x,y]$. From now,
all the groups which we consider will be finite.

Given two subgroups $H$ and $K$ of $G$ and two integers $n,m\geq1$,
we define
\begin{equation}\label{e:1} \mathrm{p}^{(n,m)}_g(H,K)=\small{\frac{|\{(x_1,\ldots, x_n, y_1, \ldots, y_m)\in H^n \times K^m \ | \ [x_1,\ldots,x_n,y_1,\ldots,y_m]=g\}|}{|H|^n \ |K|^m}}\end{equation}
as the {\it probability that a randomly chosen commutator of weight
$n+m$ of $H \times K$ is equal to a given element of $G$}. Denoting
\begin{equation}\mathcal{A}=\{(x_1,\ldots,x_n,y_1,\ldots,y_m)\in H^n \times K^m
\ | \ [x_1,\ldots,x_n,y_1,\ldots,y_m]=g\},\end{equation}
$|\mathcal{A}|=|H|^n \cdot |K|^m \cdot \mathrm{p}^{(n,m)}_g(H,K)$.
The case $n=m=1$ can be found in \cite{das-nath1} and is called
\textit{generalized commutativity degree of $G$}. For $n=m=1$ and
$H=K=G$,
\begin{equation}\label{e:2} \mathrm{p}^{(1,1)}_g(G,G)=\mathrm{p}_g(G)=\frac{|\{(x,y)\in G^2 \ | \ [x,y]=g\}|}{|G|^2}\end{equation}
is the {\it probability that the commutator of two group elements of
$G$ is equal to a given element of $G$} in \cite{ps}.

It is well known  (see  for instance \cite[Excercise 3, p. 183]{alp}) that the function $\psi(g)=|\{(x, y)\in
G\times G \ | \ [x,y]=g\}|$ is a character of $G$ and we have $\psi={\underset{\chi\in
\mathrm{Irr}(G)}\sum}\frac{|G|}{\chi(1)}\chi$, where $\mathrm{Irr}(G)$ denotes the set of all irreducible
complex characters of $G$. However, the authors exploited   this fact in \cite[Theorem 2.1]{ps}, writing
(\ref{e:2}) as
\begin{equation}\label{e:3} \mathrm{p}_g(G)= \frac{1}{|G|}{\underset{\chi\in \mathrm{Irr}(G)}\sum}\frac{\chi(g)}{\chi(1)},\end{equation}
For terminology and notations in character theory we refer to \cite{isaacs}.

Now for $g=1$,
\begin{equation}\label{e:4} \mathrm{p}^{(1,1)}_1(G,G)=\mathrm{p}_1(G)=\mathrm{d}(G)=\frac{|\{(x,y)\in G^2 \ | \ [x,y]=1\}|}{|G|^2}= \frac{|{\rm
Irr}(G)|}{|G|}\end{equation} is the $probability$ $of$ $commuting$ $pairs$ $of$ $G$ (or briefly the
$commutativity$ $degree$ of $G$), largely studied in \cite{cms, das-nath2, das-nath1, doostie-maghasedi, elr,
err, erovenko-sury, gallagher, gr, gustafson, l, rezaei1, rusin}. In particular,
\begin{equation}\label{e:5} \mathrm{p}^{(n,1)}_1(G,G)=\frac{|\{(x_1,\ldots,x_n,x_{n+1})\in G^{n+1} \ | \ [x_1,\ldots,x_n,x_{n+1}]=1\}|}{|G|^{n+1}}= \mathrm{d}^{(n)}(G),\end{equation}
is the $n$-th {\it nilpotency degree} of $G$ in \cite{cms, elr, err,
rezaei1, rezaei2} and that
\begin{equation}\label{e:6}
\mathrm{p}^{(n,1)}_1(H,G)=\frac{|\{(x_1,\ldots,x_n,y)\in H^n \times
G \ | \ [x_1,\ldots,x_n,y]=1\}|}{|H|^n \
|G|}=\mathrm{d}^{(n)}(H,G)\end{equation} is the {\it relative}
$n$-th {\it nilpotency degree} of $H$ in $G$, studied in \cite{elr,
err, rezaei1, rezaei2}. We may express (\ref{e:6})  not necessarily
with $g=1$; assuming that $H$ is normal in $G$, \cite[Equation (4)
and Theorem 4.2]{das-nath1} imply
\begin{equation}\label{e:7}\mathrm{p}^{(1,1)}_g(H,G)=\frac{|\{(x,y)\in H \times G \ | \
 [x,y]=g\}|}{|H| \ |G|}
=\frac{1}{|H||G|}{\underset{\chi \in \mathrm{Irr}(G)}\sum} \frac{|H| \langle \chi_H, \chi_H \rangle}{\chi(1)}
\chi(g),\end{equation} where $\chi_H$ denotes the restriction of $\chi$ to $H$ and $\langle , \rangle$ the usual
inner product. Our purpose is to study (\ref{e:1}), extending the previous contributions in \cite{cms,
das-nath1, elr, ps, rezaei1}. The main results of the present paper are in Section 3, in which the general
considerations of Section 2 are applied.

\section{Technical properties and some computations}
 We begin with two elementary observations on (\ref{e:1}).

 \begin{rem}\label{r:1}If $\mathcal{S}=\{[x_1,\ldots,x_n,y_1,\ldots,y_m] \ | \ x_1,\ldots,x_n \in H; y_1,\ldots,y_m \in K\}$, then  $\mathrm{p}^{(n,m)}_g (H,K)=0$ if and only if $g\not\in \mathcal{S}$. On another hand,
$\mathrm{p}^{(n,m)}_1 (H,K)=1$ if and only if
$[\underbrace{H,\ldots,H}_{n-\mathrm{times}},\underbrace{K,\ldots,K}_{m-\mathrm{times}}]=[_nH,_mK]=1$.
 \end{rem}

\begin{rem}\label{r:2}The equation (\ref{e:1}) assigns by default the map
\begin{equation}\label{e:8}\mathrm{p}^{(n,m)}_g : (x_1,\ldots,x_n,y_1,\ldots,y_m) \in H^n
\times K^m \mapsto \mathrm{p}^{(n,m)}_g(H,K)\in [0,1],\end{equation}
which is a probability measure on $H^n \times K^m$, satisfying a
series of standard properties such as being multiplicative,
symmetric and monotone.
\end{rem}

The fact that (\ref{e:8}) is multiplicative is described by the next result.

\begin{prop}\label{p:1} Let $E$ and $F$ be two groups such that $e\in E$, $f\in F$, $A,C\leq E$ and $B,D\leq F$.  Then
\[\mathrm{p}^{(n,m)}_{(e,f)} (A\times C,B\times D)=\mathrm{p}^{(n,m)}_e (A,B) \cdot \mathrm{p}^{(n,m)}_f(C,D).\]\end{prop}

\begin{proof}
It is enough to note that \[[([a_1,\ldots,a_n],[c_1,\ldots,c_n]),([b_1,\ldots,b_m],[d_1,\ldots,d_m])]
=([[a_1,\ldots,a_n],[b_1,\ldots,b_m]],[c_1,\ldots,c   _n],[d_1,\ldots,d_m]]).\]
\end{proof} Proposition \ref{p:1} is true for finitely many factors instead of only two
factors and this can be checked with easy computations. Therefore
the proof is omitted. The fact that (\ref{e:8}) is symmetric is
described by the next result.

\begin{prop}\label{p:2}With the notations of {\rm (\ref{e:1})},
$\mathrm{p}^{(n,m)}_g (H,K)=\mathrm{p}^{(n,m)}_{g^{-1}} (K,H)$.
Moreover, if $H$, or $K$, is normal in $G$, then
$\mathrm{p}^{(n,m)}_g (H,K)=\mathrm{p}^{(n,m)}_g
(K,H)=\mathrm{p}^{(n,m)}_{g^{-1}} (H,K)$.
\end{prop}

\begin{proof}
The commutator rule $[x,y]^{-1} = [y,x]$ implies the first part of
the result. Now let $H$ be normal in $G$, $n\leq m$ and
$\mathcal{B}=\{(y_1,\ldots,y_m,x_1,\ldots,x_n)\in K^m \times H^n \ |
\ [y_1,\ldots,y_m,x_1,\ldots,x_n]=g\}$. The map $\varphi :
(x_1,\ldots,x_n,y_1,\ldots,y_m) \in \mathcal{A} \mapsto
(y^{-1}_1,y^{-1}_2,\ldots,y^{-1}_n,y^{-1}_{n+1},\ldots,y^{-1}_m,y_1x_1y^{-1}_1,y_2x_2y^{-1}_2,\ldots,y_nx_ny^{-1}_n)
\in \mathcal{B}$ is  bijective and so the remaining equalities
follow. A similar argument can be applied, when the assumption $H$
is normal in $G$ is replaced by $K$ is normal in $G$.
\end{proof}

The fact that  (\ref{e:8}) is monotone is more delicate to prove,
since this is a situation in which we may find upper bounds for
(\ref{e:1}). Details are given later on. Now we will get another
expression for (\ref{e:1}). With the notations of (\ref{e:1}),
$\mathrm{Cl}_K([x_1,\ldots,x_n])$ denotes the {\it $K$-conjugacy
class of} $[x_1,\ldots,x_n] \in H$.

\begin{prop}\label{p:3}With the notations of {\rm (\ref{e:1})}, \begin{equation}\mathrm{p}_g^{(n,m)}(H,K)=\frac{1}{|H|^n \ |K|^m}  \underset{\underset{g^{-1}[x_1,\ldots,x_n]\in \mathrm{Cl}_K([x_1,\ldots,x_n])} {x_1,\ldots,x_n\in H}}\sum
|C_K([x_1,\ldots,x_n])|^m.\end{equation}
\end{prop}

\begin{proof} It is straightforward to check that \begin{equation}C_{K^m}([x_1,\ldots,x_n])=\underbrace{C_K([x_1,\ldots,x_n]) \times \ldots \times
C_K([x_1,\ldots,x_n])}_{m-\mathrm{times}}.\end{equation} In
particular, $|C_{K^m}([x_1,\ldots,x_n])|=|C_K([x_1,\ldots,x_n])|^m$.

 $\mathcal{A}=\underset{[x_1,\ldots,x_n]\in H}\bigcup \{[x_1,\ldots,x_n]\} \times T_{[x_1,\ldots,x_n]},$ where
$T_{[x_1,\ldots,x_n]}=\{(y_1,\ldots,y_m)\in K^m \ | \ [x_1,\ldots,x_n,y_1,\ldots,y_m]=g\}$. Obviously,
$T_{[x_1,\ldots,x_n]}\not=\emptyset$ if and only if $g^{-1}[x_1,\ldots,x_n]\in \mathrm{Cl}_K([x_1,\ldots,x_n])$.
Let $T_{[x_1,\ldots,x_n]}\not=\emptyset$. Then $|T_{[x_1,\ldots,x_n]}|=|C_{K^m}([x_1,\ldots,x_n])|$, because the
 map $\psi : [y_1,\ldots,y_m] \mapsto g
\overline{[y_1,\ldots,y_m]}^{^{-1}}[y_1,\ldots,y_m]$ is bijective,
where $\overline{[y_1,\ldots,y_m]}$ is a fixed element of
$T_{[x_1,\ldots,x_n]}$. We deduce that
\begin{equation}\begin{array}{lcl}|\mathcal{A}|=\sum_{[x_1,\ldots,x_n]\in H}|T_{[x_1,\ldots,x_n]}|=\underset{\underset{g^{-1}[x_1,\ldots,x_n]\in \mathrm{Cl}_K([x_1,\ldots,x_n])} {x_1,\ldots,x_n\in H}}\sum
|C_{K^m}([x_1,\ldots,x_n])|\vspace{0.3cm}\\
=\underset{\underset{g^{-1}[x_1,\ldots,x_n]\in
\mathrm{Cl}_K([x_1,\ldots,x_n])} {x_1,\ldots,x_n\in H}}\sum
|C_K([x_1,\ldots,x_n])|^m\end{array}\end{equation} and the result
follows.
\end{proof}

Special cases  of Proposition \ref{p:3} are listed below.

\begin{cor}\label{c:1}In {\rm Proposition \ref{p:3}} , if $m=1$ and $G=K$, then
\begin{equation}\mathrm{p}_g^{(n,1)}(H,G)=\frac{1}{|H|^n \ |G|}
\underset{\underset{g^{-1}[x_1,\ldots,x_n]\in
\mathrm{Cl}_G([x_1,\ldots,x_n])} {x_1,\ldots,x_n\in H}}\sum
|C_G([x_1,\ldots,x_n])|.\end{equation}
\end{cor}

\begin{cor}[See \cite{das-nath1}, Theorem 2.3]\label{c:2} In {\rm Proposition \ref{p:3}} , if $m=n=1$, then
\begin{equation}\mathrm{p}_g^{(1,1)}(H,K)=\frac{1}{|H| \ |K|}
\underset{\underset{g^{-1}x\in \mathrm{Cl}_K(x)} {x\in H}}\sum
|C_K(x)|.\end{equation} In particular, if $G=K$, then
$\mathrm{p}_g^{(1,1)}(H,G)=\frac{1}{|H| \ |G|}
\underset{\underset{g^{-1}x\in \mathrm{Cl}_G(x)} {x\in H}}\sum
|C_G(x)|$.
\end{cor}

\begin{cor}[See \cite{elr}, Proof of Lemma 4.2]\label{c:3}In {\rm Proposition \ref{p:3}} , if $m=1$ and $G=K$, then
\begin{equation}\mathrm{p}_1^{(n,1)}(H,G)=\mathrm{d}^{(n)}(H,G)=\frac{1}{|H|^n \ |G|}  \underset{x_1,\ldots,x_n\in H}\sum
|C_G([x_1,\ldots,x_n])|.\end{equation}
\end{cor}

\begin{cor}\label{c:4}In {\rm Proposition \ref{p:3}} , if $C_K([x_1,\ldots,x_n]) = 1$, then
\begin{equation}\mathrm{p}^{(n,m)}_1(H,K)=\frac{1}{|H|^n}+\frac{1}{|K|^m}-\frac{1}{|H|^n \ |K|^m}.\end{equation}
\end{cor}
\cite[Proposition 3.4]{das-nath1} follows from Corollary \ref{c:4},
when $m=n=1$.

\begin{rem}\label{r:3} Equation (\ref{e:6}) makes equivalent the study of $\mathrm{p}^{(n,1)}_1(H,G)$ and that of
$\mathrm{d}^{(n)}(H,G)$. This  is illustrated in Corollary \ref{c:3}
and noted here for the first time. Therefore there are many
information from \cite{cms, elr, err, rezaei1} and \cite{das-nath1,
das-nath2, ps} which can be connected. It is relevant to point out
that these concepts were treated independently and with different
methods in the last years.
\end{rem}

Let $\chi$ be a character of $G$ and
$\theta$ be a character of $H\leq G$. The \textit{Frobenius
Reciprocity Law} \cite[Lemma 5.2]{isaacs} gives a link between the
restriction $\chi_H$ of $\chi$ to $H$ and the induced character
$\theta^G$ of $\theta$. Therefore $\langle
\chi,\theta^G\rangle_G=\langle \chi_{H},\theta\rangle_{H}.$Write
this number as
$e_{(\chi,\theta)}=\langle \chi,\theta^G\rangle_G=\langle \chi_H,\theta\rangle_H.$
If $e_{(\chi, \theta)}=0$, then $\theta$ does not appear in $\chi_H$ and so $\chi$ does not appear in
$\theta^G$. Recall from \cite{isaacs} that, if $e_{(\chi,\theta)}\neq 0$, then $\chi$ \textit{covers} $\theta$
(or also $\theta$ \textit{belongs to the constituents of} $\chi_H$). In particular, if $\theta= \chi_{H}$, then
$e_{(\chi,\chi_H)}=\langle \chi,(\chi_H)^G\rangle_{G}=\langle \chi_H,\chi_H\rangle_H.$
From a classic relation (see \cite[Lemma 2.29]{isaacs}),
$e_{(\chi,\chi_H)}=\langle \chi,(\chi_H)^G\rangle_G=\langle \chi_H,\chi_H\rangle_H\leq |G:H| \ \langle \chi,\chi
\rangle_G = |G:H| e_{(\chi,\chi)}$
 and the equality holds if and only if $\chi(x)=0$ for all $x\in G-H$. In particular, if $\chi \in \mathrm{Irr}(G)$, then
$\langle \chi_H,\chi_H\rangle_H= |G:H| \ \mathrm{if \ and \ only \ if} \ \chi(x)=0,$ for all $x\in G-H.$
Therefore the following result is straightforward.

\begin{cor}With the notations of {\rm (\ref{e:1})},
$\mathrm{p}_{g}^{(1,1)}(H,G)\leq |G:H| \ \mathrm{p}_1(G)$ and the equality holds if and only if all the
characters vanish on $G-H$.
\end{cor}



At this point, \cite[Theorem 4.2]{das-nath1} becomes
\begin{equation}\label{e:21}
\zeta(g)=|H| \ \underset{\chi \in \mathrm{Irr}(G)}\sum \frac{e_{(\chi_H,\chi_H)}}{\chi(1)} \cdotp
\chi(g)=|\{(x,y)\in H \times G \ | \
 [x,y]=g\}|= \underset{\underset{g^{-1}x \in \mathrm{Cl}_G(x)} {x\in
H}}\sum  |C_G(x)|,
\end{equation}
where $\zeta(g)$ is the number of solutions $(x,y)\in H\times G$ of the equation $[x,y]=g$. Note that
(\ref{e:21}) and \cite[Excercise 3, p. 183]{alp} give a short argument  to prove that $\zeta(g)$ is a character
of $G$ with respect to the argument in \cite[Corollary 4.3]{das-nath1}. The equation (\ref{e:7}) becomes
\begin{equation}\label{e:22}
\mathrm{p}_{g}^{(1,1)}(H,G)=\frac{\zeta(g)}{|H| \ |G|}.
\end{equation}


For the general case that $n>1$, $m>1$ and $G=K$,
\begin{equation}\label{e:23}
\mathrm{p}_{g}^{(n,m)}(H, G)=\frac{\zeta^{(n,m)}(g)}{|G|^m}= \frac{1}{|G|^m} \Big(
\underset{\underset{g^{-1}[x_1,\ldots,x_n]\in \mathrm{Cl}_G([x_1,\ldots,x_n])} {x_1,\ldots,x_n\in H}}\sum
|C_G([x_1,\ldots,x_n])|^m\Big), \end{equation} where
\begin{equation}\label{e:extra}
\zeta^{(n,m)}(g)= \underset{\underset{g^{-1}[x_1,\ldots,x_n]\in \mathrm{Cl}_G([x_1,\ldots,x_n])}
{x_1,\ldots,x_n\in H}}\sum |C_G([x_1,\ldots,x_n])|^m
\end{equation}
is the number of solutions $(x_1,\ldots,x_n,y_1,\ldots,y_m)\in H^n\times G^m$ of
$[x_1,\ldots,x_n,y_1,\ldots,y_m]=g$.


\begin{rem}\label{t:1}There are many evidences from the computations that $\zeta^{(n,m)}(g)$ is a character of $G$.
\end{rem}




Now we may prove upper bounds for (\ref{e:1}) and find that
(\ref{e:8}) is monotone.

\begin{prop}\label{p:4} With the notations of
{\rm (\ref{e:1})}, if $H \leq K$, then
$\mathrm{p}^{(n,m)}_g(H,G)\geq \mathrm{p}^{(n,m)}_g(K,G).$ The
equality holds if and only if $\mathrm{Cl}_H(x)=\mathrm{Cl}_K(x)$
for all $x\in G$.
\end{prop}
\begin{proof}
We note that $\frac{1}{|K|}\leq \frac{1}{|H|}$ and then
$\frac{1}{|K|^n}\leq \frac{1}{|H|^n}$. By Proposition \ref{p:3},
\begin{equation}\label{e:24}\begin{array}{lcl}
|G|^m \cdot \mathrm{p}_{g}^{(n,m)}(K,G)=\frac{1}{|K|^n}
\underset{\underset{g^{-1}[x_1,\ldots,x_n]\in
\mathrm{Cl}_G([x_1,\ldots,x_n])} {x_1,\ldots,x_n\in K}}\sum
|C_G([x_1,\ldots,x_n])|\vspace{0.3cm}\\
\leq \frac{1}{|H|^n} \underset{\underset{g^{-1}[x_1,\ldots,x_n]\in
\mathrm{Cl}_G([x_1,\ldots,x_n])} {x_1,\ldots,x_n\in K}}\sum
|C_G([x_1,\ldots,x_n])|\end{array}\end{equation} in particular the
last relation is true for $x_1,\ldots,x_n\in H\leq K$ and continuing
\begin{equation}\label{e:25}=\frac{1}{|H|^n}
\underset{\underset{g^{-1}[x_1,\ldots,x_n]\in \mathrm{Cl}_G([x_1,\ldots,x_n])} {x_1,\ldots,x_n\in H}}\sum
|C_G([x_1,\ldots,x_n])|=|G|^m \cdot p_{g}^{(n, m)}(H, G).
\end{equation}
The rest of the proof is clear.
\end{proof}

The next result shows an upper bound, which generalizes
\cite[Theorem 4.6]{elr}.

\begin{prop}\label{p:5} With the notations of {\rm (\ref{e:1})}, if $N$ is a normal subgroup of $G$ such that $H\leq N$, then
$\mathrm{p}^{(n,m)}_g(H,G)\leq \mathrm{p}^{(n,m)}_g
\Big(\frac{H}{N},\frac{G}{N}\Big).$ Moreover, if $N \cap [_nH,_mG] =
1$, then the equality holds.
\end{prop}
\begin{proof}We have
\[\begin{array}{lcl}
|H|^n \ |G|^m \ \mathrm{p}^{(n,m)}_g(H,G)=|\mathcal{A}|\vspace{0.3cm}\\
=|\{(x_1,\ldots,x_n,y_1,\ldots,y_m) \in H^n \times G^m \ | \
[x_1,\ldots,x_n,y_1,\ldots,y_m] \cdot
g^{-1}=1 \}|\vspace{0.3cm}\\
=|\{(x_1,\ldots,x_n,y_1,\ldots,y_m) \in H^n \times G^m \ | \
[x_1,\ldots,x_n,y_1,\ldots,y_m,g^{-1}]=1 \}|\vspace{0.3cm}\\
 =\sum_{x_1 \in H} \ldots \sum_{x_n\in
 H}\sum_{y_1 \in G} \ldots \sum_{y_m\in
 G}|C_G([x_1,\ldots,x_n,y_1,\ldots,y_m])|\vspace{0.3cm}\\
 =\sum_{x_1 \in H} \ldots \sum_{x_n\in
 H}\sum_{y_1 \in G} \ldots \sum_{y_m\in
 G}\frac{|C_G([x_1,\ldots,x_n,y_1,\ldots,y_m])N| \cdot |C_N([x_1,\ldots,x_n,y_1,\ldots,y_m])|}{|N|}\vspace{0.3cm}\\
\leq \sum_{x_1 \in H} \ldots \sum_{x_n\in
 H}\sum_{y_1 \in G} \ldots \sum_{y_m\in
 G}|C_{G/N}([x_1N,\ldots,x_nN,y_1N,\ldots,y_mN])| \vspace{0.3cm}\\
 \cdot |C_N([x_1,\ldots,x_n,y_1,\ldots,y_m])|\vspace{0.3cm}\\
= \sum_{S_1 \in H/N} \sum_{x_1 \in S_1}\ldots \sum_{S_n\in
 H/N}\sum_{x_n \in S_n}\sum_{T_1 \in G/N} \sum_{y_1 \in T_1} \ldots \sum_{T_m\in
 G/N}\sum_{y_m \in T_m}\vspace{0.3cm}\\
 |C_{G/N}([S_1,\ldots,S_n,T_1,\ldots,T_m])| \cdot |C_N([x_1,\ldots,y_m])|\vspace{0.3cm}\\
= \Big(\sum_{S_1 \in H/N} \ldots \sum_{S_n\in
 H/N}\sum_{T_1 \in G/N} \ldots \sum_{T_m\in
 G/N}|C_{G/N}([S_1,\ldots,S_n,T_1,\ldots,T_m])|\Big)\vspace{0.3cm}\\
 \cdot \Big(\sum_{x_1 \in S_1}\ldots \sum_{x_n \in S_n}\sum_{y_1 \in T_1} \ldots \sum_{y_m \in T_m}
 |C_N([x_1,\ldots,y_m])|\Big)\vspace{0.3cm}\\
\leq |N|^{n+m} \sum_{S_1 \in H/N} \ldots \sum_{S_n \in H/N}\sum_{T_1
\in G/N} \ldots\sum_{T_m\in
 G/N}\vspace{0.3cm}\\
 |C_{G/N}([S_1,\ldots,S_n,T_1,\ldots,T_m])|\vspace{0.3cm}\\
= \Big|\frac{H}{N}\Big|^n \ \Big|\frac{G}{N}\Big|^m \
p^{(n,m)}_g\Big(\frac{H}{N},\frac{G}{N}\Big) \ |N|^{n+m}=|H|^n \
|G|^m \ \mathrm{p}^{(n,m)}_g\Big(\frac{H}{N},\frac{G}{N}\Big).
\end{array}\]
The condition of equality in the above relations is satisfied
exactly when $N \cap [_nH,_mG] = 1$. The result follows.\end{proof}

\begin{cor} A special case of {\rm Proposition \ref{p:5}}  is $\mathrm{p}_g(G)\leq \mathrm{p}_g(G/N)$.
\end{cor}

\begin{cor}[See \cite{elr}, Theorem 4.6] In {\rm Proposition \ref{p:5}} , if $m=1$ and $g=1$, then $\mathrm{d}^{(n)}(H,G)\leq \mathrm{d}^{(n)}(H/N,G/N)$.
\end{cor}

\section{Some upper and lower bounds}

A relation among (\ref{e:1})--(\ref{e:7}) is described below.

\begin{thm}\label{t:2} With the notations of {\rm (\ref{e:1})},
$\mathrm{p}^{(n,m)}_g(G,G)\leq \mathrm{p}^{(n,m)}_g(H,K)\leq
\mathrm{p}^{(n,m)}_1(H,K)\leq \mathrm{p}^{(n,m)}_1(H,G) \leq
\mathrm{p}^{(n,m)}_1(H,H).$
\end{thm}
\begin{proof}
From  Proposition \ref{p:4},  $\mathrm{p}^{(n,m)}_g(G,G)\leq
\mathrm{p}^{(n,m)}_g(G,H)$. From Proposition \ref{p:3},
\begin{equation}
\mathrm{p}^{(n,m)}_g(H,K)=\frac{1}{|H|^n \ |K|^m}
\underset{\underset{g^{-1}[x_1,\ldots,x_n]\in
\mathrm{Cl}_K([x_1,\ldots,x_n])} {x_1,\ldots,x_n\in H}}\sum
|C_K([x_1,\ldots,x_n])|^m
\end{equation} and for $g=1$ we get \begin{equation}\leq \frac{1}{|H|^n \ |K|^m}
\underset{x_1,\ldots,x_n\in H}\sum
|C_K([x_1,\ldots,x_n])|^m=\mathrm{p}^{(n,m)}_1(H,K),
\end{equation}where in the last passage still Proposition \ref{p:3} is used.
From $C_K([x_1,\ldots,x_n])\subseteq C_G([x_1,\ldots,x_n])$, we
deduce
\begin{equation}\leq
\underset{x_1,\ldots,x_n\in H}\sum |C_G([x_1,\ldots,x_n])|^m=
\mathrm{p}^{(n,m)}_1(H,G).
\end{equation}
Applying  Proposition \ref{p:2},
$\mathrm{p}^{(n,m)}_1(H,G)=\mathrm{p}^{(n,m)}_1(G,H)$ and so
$\mathrm{p}^{(n,m)}_1(G,H)\leq \mathrm{p}^{(n,m)}_1(H,H)$ by
Proposition \ref{p:4}.
\end{proof}

\begin{cor}\label{c:5} With the notations of {\rm (\ref{e:1})}, if $Z(G)=1$, then
$\mathrm{p}^{(n,1)}_g(H,K)\leq \frac{2^n-1}{2^n}.$\end{cor}

\begin{proof} It follows from Theorem \ref{t:2} and \cite[Theorem 5.3]{elr}.
\end{proof}

Another significant restriction is the following.

\begin{thm}\label{t:3}With the notations of {\rm (\ref{e:1})}, let $p$ be the smallest prime divisor of $|G|$. Then
\begin{itemize}
\item[(i)]$\mathrm{p}^{(n,m)}_g(H,K)\leq  \frac{2p^n+p-2}{p^{m+n}}$;
\item[(ii)]$\mathrm{p}^{(n,m)}_g(H,K)\geq \frac{(1-p)|Y_{H^n}|+p|H^n|}{|H^n| \ |K^m|} - \frac{(|K|+p) |C_H(K)|^n}{|H^n| \ |K^m|}$;
\end{itemize}
where $Y_{H^n}=\{[x_1,\ldots,x_n]\in H^n \ | \ C_K([x_1,\ldots,x_n])=1\}$.
\end{thm}

\begin{proof} If $[_nH,_mK]=1$, then $C_{H^n}(K^m)=H^n$ and $Y_{H^n}$ equals $H^n$ or an empty set according as $K^m$ is trivial or nontrivial.
Assume that $[_nH,_mK]\not=1$. Then $Y_{H^n}\cap C_{H^n}(K^m)=Y_{H^n}\cap (C_H(K^m)\times \ldots \times
C_H(K^m))=Y_{H^n}\cap (C_H(K)\times C_H(K)\times \ldots \times C_H(K))=Y_{H^n}\cap (C_H(K))^{nm}\not=\emptyset$
and \begin{equation}\begin{array}{lcl} \underset{x_1,\ldots,x_n\in H}\sum |C_{K^m}([x_1,\ldots,x_n])|=
\underset{x_1,\ldots,x_n\in H}\sum |C_K([x_1,\ldots,x_n])|^m\vspace{0.3cm}\\
 =\underset{x_1,\ldots,x_n\in Y_{H^n}}\sum
|C_K([x_1,\ldots,x_n])|^m +  \underset{x_1,\ldots,x_n\in
C_{H^n}(K)}\sum |C_K([x_1,\ldots,x_n])|^m \vspace{0.3cm}\\
+\underset{x_1,\ldots,x_n\in H^n-(Y_{H^n}\cup C_{H^n}(K))}\sum
|C_K([x_1,\ldots,x_n])|^m \vspace{0.3cm}\\
=|Y_{H^n}|+ |K| \ |C_H(K)|^n+\underset{x_1,\ldots,x_n\in H^n-(Y_{H^n}\cup C_{H^n}(K))}\sum
|C_K([x_1,\ldots,x_n])|^m.\end{array}\end{equation} Since $p^m\leq |C_K([x_1,\ldots,x_n])|^m\leq
\frac{|K^m|}{p^m}$, $|Y_{H^n}|\leq |H^n|$ and $p^n\leq|C_H(K)|^n\leq \frac{|H^n|}{p^n}$, \begin{equation}\leq
|Y_{H^n}|+|K| \ |C_H(K)|^n+(|H^n|-(|Y_{H^n}|+ |C_H(K)|^n) \cdot \frac{|K^m|}{p^m}\end{equation} and then
\begin{equation}\label{e:bound}\begin{array}{lcl}\mathrm{p}^{(n,m)}_g(H,K)\leq \frac{|Y_{H^n}|}{|H^n| \ |K^m|} + \frac{|K| \
|C_H(K)|^n}{|H^n| \ |K^m|}+ \frac{1}{p^m} - \frac{|Y_{H^n}|}{p^m \ |H^n|} -
\frac{|C_H(K)|^n}{p^m \ |H^n|}\vspace{0.3cm}\\
\leq  \frac{1}{p^m} + \frac{1}{p^{m+n-1}} + \frac{1}{p^m}-\frac{1}{p^{m+n}}-\frac{1}{p^{m+n}} =
\frac{2p^n+p-2}{p^{m+n}}.\end{array}\end{equation} Hence (i) follows. On another hand, we may continue in the
other direction \begin{equation}\geq |Y_{H^n}|+|K| \ |C_H(K)|^n+ p \ (|H^n|-(|Y_{H^n}|+
|C_H(K)|^n)\end{equation} and then \begin{equation}\mathrm{p}^{(n,m)}_g(H,K)\geq
 \frac{(1-p)|Y_{H^n}|}{|H^n| \ |K^m|} + \frac{p}{|K^m|} - \frac{(|K|+p) |C_H(K)|^n}{|H^n| \ |K^m|}.\end{equation} Then (ii)
follows.
\end{proof}

The bound in Theorem \ref{t:3} (i) is a little bit different from
the bound in \cite[Corollary 3.9]{das-nath1}, where it is proved
that $\mathrm{p}^{(1,1)}_g(H,K)\leq \frac{2p-1}{p^2}$ and in
particular $\mathrm{p}^{(1,1)}_g(H,K)\leq \frac{3}{4}$. We conclude
the following structural restriction.

\begin{cor}\label{c:6}
In {\rm Theorem} \ref{t:3}, if $\mathrm{p}^{(n,m)}_g(H,K)= \frac{2p^n+p-2}{p^{m+n}}$, then
\begin{equation}|H:C_H(K)|\leq \Big(\frac{p^{n+1}-p^3-\frac{p^2}{2}+p}{2p^2+p-2}\Big)^{\frac{1}{n}}.\end{equation}
\end{cor}

\begin{proof} Looking at (\ref{e:bound}) and the proof of Theorem \ref{t:3} (i), we deduce
\begin{equation}\begin{array}{lcl}
\frac{2p^n+p-2}{p^{m+n}}\leq \frac{|Y_{H^n}|}{|H^n| \ |K^m|} + \frac{|K| |C_H(K)|^n}{|H^n| \ |K^m|} +
\frac{1}{p^m}\leq \frac{1}{p^m} + \frac{1}{p^{m-1}} \ \Big|\frac{C_H(K)}{H}\Big|^n + \frac{1}{p^m}\vspace{0.3cm}\\
=\frac{1}{p^{m-1}} \Big( \frac{2}{p}+|\frac{C_H(K)}{H}|^n\Big)
\end{array}\end{equation}
and then $\frac{2p^n+p-2}{p^{n+1}}\leq \frac{2}{p}+\Big|\frac{C_H(K)}{H}\Big|^n$. We conclude that
$\frac{p^{n+1}}{2p^n+p-2}\geq \frac{p}{2}+\Big|\frac{H}{C_H(K)}\Big|^n$ and so
\begin{equation}\frac{p^{n+1}}{2p^n+p-2}- \frac{p}{2} =\frac{p^{n+1}-p^3-\frac{p^2}{2}+p}{2p^2+p-2}\geq
\Big|\frac{H}{C_H(K)}\Big|^n.\end{equation} The result follows, once we extract the $n$-th root.
\end{proof}

\section*{Acknowledgement}
The second author is grateful to the colleagues of the Ferdowsi University of Mashhad  for some helpful comments
in the period  in which the present work has been written.


\begin{thebibliography}{20}

\bibitem {alp} J. Alperin and B. Bell, \textit{Groups and Representations}, Springer, 1995, New York.

\bibitem{cms}  K. Chiti, M. R. R. Moghaddam and A. R. Salemkar, $n$--isoclinism classes and $n$--nilpotency degree of finite groups, {\it Algebra Colloq.} {\bf 12} (2005), 225--261.


\bibitem{das-nath2} A. K. Das and R. K. Nath, On solutions of a class of equations in a finite group, \textit{Commun. Algebra} \textbf{37} (2009), 3904--3911.

\bibitem{das-nath1}
A. K. Das and R. K. Nath, On the generalized relative commutative
degree of a finite group, \textit{Int. Electr. J. Algebra}
\textbf{7} (2010), 140--151.

\bibitem{doostie-maghasedi}H. Doostie and  M. Maghasedi, Certain classes of groups with commutativity degree $d(G)<
1/2$, \textit{Ars Combinatoria} \textbf{89} (2008), 263--270.


\bibitem{et}   P. Erd\H{o}s and P. Tur\'an, On some problems of statistical group theory, {\it Acta Math. Acad. Sci. Hung.} {\bf 19} (1968), 413--435.

\bibitem{elr} A. Erfanian,  P. Lescot and  R. Rezaei, On the relative commutativity degree of a subgroup of a finite group, \textit{Comm. Algebra} {\bf{35}} (2007), 4183--4197.

\bibitem{er} A. Erfanian and R. Rezaei, On the commutativity degree of compact groups, \textit{Arch. Math. (Basel)} \textbf{93} (2009), 345--356.

\bibitem{err} A. Erfanian, R. Rezaei and F. G. Russo, Relative $n$-isoclinism classes and relative $n$-th nilpotency degree of finite groups, e-print, Cornell University, 2010, arXiv:0003310 [math.GR].





\bibitem{erovenko-sury}I.V. Erovenko and B. Sury, Commutativity degree of wreath products of finite abelian groups, \textit{Bull. Austral. Math. Soc.} \textbf{77} (2008), 31--36.

\bibitem{gallagher} P. X. Gallagher,  The number of conjugacy classes in a finite group, \textit{Math. Z.} {\bf 118} (1970), 175--179.

\bibitem{gr} R. M. Guralnick and G. R. Robinson, On the commuting probability in finite groups, \textit{J. Algebra} \textbf{300} (2006), 509--528.

\bibitem{gustafson}W. H. Gustafson,  What is the probability that two groups elements commute?  \textit{Amer. Math. Monthly} {\bf 80} (1973), 1031--1304.

\bibitem{isaacs}I. M. Isaacs, {\it Character Theory of Finite Groups}, Dover Publ., New York, 1994.

\bibitem{l} P. Lescot, Isoclinism classes and commutativity degrees of finite groups, \textit{J. Algebra} {\bf{177}} (1987), 847--869.




\bibitem{ps} M. R. Pournaki and R. Sobhani, Probability that the commutator of two group elements is equal to a given element, \textit{J. Pure Appl. Algebra} \textbf{212} (2008), 727--734.

\bibitem{rezaei1} R. Rezaei and F. G. Russo,   $n$-th relative nilpotency degree and relative $n$-isoclinism classes, e-print, Cornell University, 2010, arXiv:1003.2297v1 [math.GR].

\bibitem{rezaei2}R. Rezaei and F. G. Russo, Bounds for the relative $n$-th nilpotency degree in compact groups, e-print, Cornell University, 2009,  arXiv:0910.4716v1 [math.GR].

\bibitem{rusin} D. J. Rusin, What is the probability that two elements of a finite group commute?, \textit{Pacific J. Math.} \textbf{82} (1979), 237--247.



\end{thebibliography}
\end{document}